\documentclass{elsart}
\usepackage{graphics}
\usepackage{epsfig}
\usepackage{amssymb}
\begin{document}
\begin{frontmatter}
\title{Application of the HLSVD technique
to the filtering of X-ray diffraction data}
\author{M. Ladisa\corauthref{cor2}}
\corauth[cor2]{Corresponding Author: Massimo Ladisa
(m.ladisa@ic.cnr.it)}
\ead{m.ladisa@ic.cnr.it} \ead[url]{www.ic.cnr.it}
\address{Istituto di Cristallografia (IC-CNR),
         Via Amendola 122/O,
         70126 Bari, Italy}
\author{A. Lamura}
\ead{a.lamura@ba.iac.cnr.it} \ead[url]{www.iac.cnr.it}
\address{Istituto Applicazioni Calcolo (IAC-CNR), Sezione di Bari,\\
          Via Amendola 122/D,
          70126 Bari, Italy}
\author{T. Laudadio}
\ead{Laudadio@esat.kuleuven.ac.be} \ead[url]{www.kuleuven.ac.be}
\address{Katholieke Universiteit Leuven,
         Department of Electrical Engineering, Division ESAT-SCD (SISTA),
         Kasteelpark Arenberg 10,
         3001 Leuven-Heverlee, Belgium}
\author{G. Nico}
\ead{g.nico@ba.iac.cnr.it} \ead[url]{www.iac.cnr.it}
\address{Istituto Applicazioni Calcolo (IAC-CNR), Sezione di Bari,\\
          Via Amendola 122/D,
          70126 Bari, Italy}

\begin{abstract}
A filter based on the Hankel Lanczos Singular Value Decomposition
(HLSVD) technique is presented and applied for the first time to
X-ray diffraction (XRD) data. Synthetic and real powder XRD
intensity profiles of nanocrystals are used to study the filter
performances with different noise levels.
Results show the robustness of the HLSVD filter and its capability to extract
easily and efficiently the useful crystallographic
information. These characteristics make the filter an interesting
and user-friendly tool for processing XRD data.
\end{abstract}


\end{frontmatter}

\section{Introduction}

In many applications of X-ray diffraction (XRD) techniques to the
study of crystal properties, a key step in the data processing
chain is an effective and adaptive noise filtering
\cite{mierzwa97,hieke99,schmidt02,rajagopal04}. A correct noise
removal can facilitate the separation of the useful
crystallographic information from the background signal, and the
estimation of crystal structure and domain size. Important issues of
XRD data filtering are performances in noise suppression,
capability to preserve the peak position, computational cost
and, finally, the possibility of being used as a blackbox tool.
Different digital filters have been applied to XRD data, in
spatial and frequency domains. Simple procedures are based on
polynomial filtering (and fitting) in the spatial domain
\cite{mierzwa97}. A standard practice when working in frequency
domain is to use Fourier smoothing. It consists in removing the
high-frequency components of the spectrum \cite{aubanel85}. Since
the truncation of high-frequency components can be problematic in
the case of high level noise, a different approach based on the
Wiener Fourier (WF) filter has been proposed to clean XRD data
\cite{wooff86}. A different approach, which makes use of the singular
value decomposition (SVD), has been successfully applied to
time-resolved XRD data to reduce noise level
\cite{schmidt02,rajagopal04}.

In this work we describe an application of the Hankel Lanczos
Singular Value Decomposition (HLSVD) algorithm to filter XRD
intensity data. The proposed filter is based on a subspace-based 
parameter estimation method, called Hankel Singular Value
Decomposition (HSVD) \cite{hsvd}, which is currently applied to
Nuclear Magnetic Resonance spectroscopy data for solvent
suppression \cite{laudadio02}. The HSVD method computes a
``signal'' subspace and a ``noise'' subspace by means of the SVD
of the Hankel matrix $H$, whose entries are the noisy signal data points.
Its computationally most intensive part consists of the
computation of the SVD of the matrix $H.$ Recently, several
improved versions of the algorithm have been developed in order to
reduce the needed computational time \cite{laudadio02}. 
In this paper, we choose to
apply the HSVD method based on the Lanczos algorithm with partial
reorthogonalization (HLSVD--PRO), which is proved to be the most
accurate and efficient version available in the literature. A
comparison in terms of numerical reliability and computational
efficiency of HSVD with its Lanczos-based variants can be found in
Ref.~\cite{laudadio02}.

A criterion is presented to facilitate the separation of noise
from the useful crystallographic signal. It enables the design of
a blackbox filter to be used in the processing of XRD data. 
Here, the filter is applied to nanocrystalline XRD data. 
Nanocrystals are characterized  by chemical and physical
properties different from those of the bulk \cite{wales96}. At a
scale of a few nanometers, metals can crystallize in a structure
different from that of bulk. Nowadays, different branches of science 
and engineering are benefiting from the properties of nanocrystalline
materials \cite{siegel}. In particular, recent XRD experiments have shown 
that intensities, measured as a function of the scattering angle, could
be useful to extract structural and domain size information about
nanocrystalline materials. Synthetic XRD datasets are generated
by computing the X-ray scattered intensity from nanocrystalline
samples of different size and properties by using an analytic
expression (see eq.~(\ref{eq:debye}) in Section 3). 
Synthetic datasets are processed and filter performance is studied 
when considering different levels of noise. 
Numerical tests on real XRD data of Au nanocrystalline samples of 
different size and properties show the
robustness of the proposed filter and its capability to extract easily and
efficiently the useful crystallographic information.
These characteristics make this filter an interesting and
user--friendly tool for the interactive processing of XRD data.

The paper has the following structure. Section 2 is devoted to the
theoretical aspects of the proposed approach. The dataset used to
study the filter properties is described in Section 3. Numerical
results are reported in Section 4. Finally, some conclusions are
drawn in Section 5.

\section{The subspace-based parameter estimation method HSVD}
\label{sec:HSVD}

Let us denote with $I_n$ the samples of the diffracted intensity
signal collected at angles $\vartheta_n$, $n = 0, \ldots, N-1$.
They are modelled as the sum of $K$ exponentially damped complex
sinusoids
\begin{equation}
I_n = I^0_n + e_n =\sum_{k=1}^K a_k \exp \left ( i \varphi_k
\right ) \exp \left [ \left (  -d_k + i 2 \pi f_k \right )
\vartheta_n \right ] + e_n
                          \;,
\label{eq:model}
\end{equation}
where $I_n$ and $I^0_n$, respectively, represent the measured and 
modelled intensities  at the {\it n-th} scattering angle $\vartheta_n = n
\Delta \vartheta + \vartheta_0$, with $\Delta \vartheta$ the
sampling angular interval and $\vartheta_0$ the initial scattering
angular position, $a_k$ is the amplitude, $\varphi_k$ the phase,
$d_k$ the damping factor and $f_k$ the frequency of the {\it k-th}
sinusoid, $k = 1, \ldots, K$, with $K$ the number of damped
sinusoids, and $e_n$ is complex white noise with a circular
Gaussian distribution. It is worth noting that the value of $K$ increases
or decreases by 2 in order to guarantee that the modelled intensity is real.

The $N$ data points defined in (\ref{eq:model}) are arranged into
a Hankel matrix $H$ of dimensions $L \times M$, with $L+M = N+1$ and
$L\simeq N/2$
\begin{equation}
H_{L \times M} = \left [
\begin{array}{cccccc}
$$I_0$$       &  $$I_1$$        & $$\ldots$$  &  & $$\ldots$$ & $$I_{M-1}$$   \\
$$I_1$$       &  $$I_2$$        & $$\ldots$$  &  & $$\ldots$$ & $$I_{M}$$     \\
$$\vdots$$    &  $$\vdots$$     & $$\vdots$$  &  & $$\vdots$$ & $$\vdots$$    \\
$$I_{L-1}$$   &  $$I_{L-2}$$    & $$\ldots$$  &  & $$\ldots$$ & $$I_{N-1}$$   \\
\end{array}
\right ]_{L \times M}. \label{eq:matrix}
\end{equation}
The SVD of the Hankel matrix  is computed as
\begin{equation}
H_{L \times M} = U_{L \times L} \Sigma_{L \times M} V^H_{M \times
M}, \label{eq:hswd}
\end{equation}
where $\Sigma = {\rm diag}\{\lambda_1, \lambda_2, \ldots,
\lambda_r\}$, $\lambda_1 \geq \lambda_2 \geq \ldots \lambda_r \geq 0$, $r
= \min(L, M)$, $U$ and $V$ are orthogonal matrices and the
superscript H denotes the Hermitian conjugate. The SVD is computed
by using the Lanczos bidiagonalization algorithm with partial
reorthogonalization \cite{simon}. 
This algorithm computes two matrix-vector products at each step.
Exploiting the structure of the matrix (\ref{eq:matrix}) by using the FFT,
the latter computation requires $O((L+M) log_2 (L+M))$ rather than $O(LM)$.

In order to obtain the ``signal'' subspace, the matrix $H$ is
truncated to a matrix $H_{K}$ of rank $K$
\begin{equation}
H_K = U_K \Sigma_K V^H_K \label{eq:hswd2},
\end{equation}
where $U_K$, $V_K$, and $\Sigma_K$ are defined by taking the first
$K$ columns of $U$ and $V$, and the $K \times K$ upper-left matrix
of $\Sigma$, respectively. As subsequent step, the least--squares
solution of the following over-determined  set of equations is
computed
\begin{equation}
V^{(top)}_K E^H \simeq V^{(bottom)}_K \label{eq:ov_eq}
\end{equation}
where $V^{(bottom)}_K$ and $V^{(top)}_K$ are derived from $V_K$ by
deleting its first and last row, respectively.
The $K$ eigenvalues ${\hat z}_k$ of the matrix $E$ are used to
estimate the frequencies ${\hat f}_k$ and the damping factors ${\hat
d}_k$ of the model damped sinusoids from the relationship
\begin{equation}
{\hat z}_k = \exp \left [ \left ( - {\hat d}_k + i 2 \pi
{\hat f}_k \right ) \Delta \vartheta \right ], \label{eq:eigen}
\end{equation}
with $k = 1, \ldots, K$. The values so obtained are inserted into the
model eq.~(\ref{eq:model}) which yields the set of equations
\begin{equation}
I_n \simeq \sum_{k=1}^K a_k \exp \left ( i \varphi_k \right )
\exp \left [ \left (  -{\hat d}_k + i 2 \pi {\hat f}_k \right
) \vartheta_n \right ] + e_n
                          \;,
\label{eq:model_eigen}
\end{equation}
with $n = 0, 1, \ldots, N-1$. The least-squares solution of
(\ref{eq:model_eigen}) provides the amplitude ${\hat a}_k$ and
phase ${\hat \varphi}_k$ estimates of the model sinusoids.

\section{Dataset}

The HLSVD--PRO filter was applied to synthetic as well as
real XRD data. In this section, the generation of XRD intensity
profiles and the experimental setup for the acquisition of real
data are described. Both synthetic and real XRD data refer to Au
nanocrystalline samples.
Nanocrystals are made of clusters of three different structure
types: cuboctahedral ${\mathcal C}$, icosahedral ${\mathcal I}$,
and decahedral ${\mathcal D}$. For each fixed structure type
${\mathcal X}$ (${\mathcal X}={\mathcal C}, {\mathcal I},
{\mathcal D}$) the size $n$ of clusters follows a log-normal
distribution
\begin{equation}
f_{\mathcal X} (n) = \frac{\exp \left ( - s_{\mathcal X} / 2
\right )}
                          {\sqrt{2 \pi \xi_{\mathcal X} s_{\mathcal X}}}
\exp{\left [ -\frac{(\log n - \log \xi_{\mathcal X})^2}
            {2 s^2_{\mathcal X}} \right ]}
                          \;,
\label{eq:nano_size_distr}
\end{equation}
with mode $\xi_{\mathcal X}$ and logarithmic width $s_{\mathcal
X}$. Structural distances for the different structure types
${\mathcal X}$ are generally studied independently of the actual
nanomaterial.
The nearest distance between atoms in the crystals is
chosen as a reference length and arbitrarily set to $1/\sqrt{2}$,
a constant in various structures ${\mathcal X}$ and for all sizes
$n$ of the clusters.
Actual distances in nanoclusters are then recovered by applying a
correction factor $a_{\mathcal X}(n)$ for strain, supposed to be
uniform and isotropic. A convenient description of the strain
factor as a function of the structure type and cluster size is
\begin{equation}
a_{\mathcal X} (n) = \Omega_{\mathcal X} +
                     \left ( \Xi_{\mathcal X} - \Omega_{\mathcal X}
                     \right ) \times
\frac{\pi + 2 \, {\rm atan} \left (
                         \frac{n^0_{\mathcal X}-n}{w_{\mathcal X}}
                         \right )}
                     {\pi + 2 \, {\rm atan} \left (
                         \frac{n^0_{\mathcal X}-1}{w_{\mathcal X}}
                         \right )}
                          \;,
\label{eq:strain_factor}
\end{equation}
given in terms of the four parameters $[n^0_{\mathcal X},
\Omega_{\mathcal X}, \Xi_{\mathcal X}, w_{\mathcal X}]$.
Intensities scattered by nanoclusters with size $n$ and structure
type ${\mathcal X}$ are computed by using the diffractional model
based on the Debye function method \cite{nano:rietveld,nano:cervellino}:
\begin{equation}
I_{{\mathcal X},n} (q) = A
                     \left \{
                    N_{\mathcal X} (n) + \sum_{i \neq j}^{N_{\mathcal X} (n)}
\frac{sin[2 \pi q u^{{\mathcal X},n}_{i,j} a_{\mathcal X} (n)]} {2
\pi q u^{{\mathcal X},n}_{i,j} a_{\mathcal X} (n)}
                      \right \}
                          \;,
\label{eq:debye}
\end{equation}
where $I_0$ is the incident X-ray intensity, $T(q')$ the
Debye-Waller factor, $f(q)$ the atomic form factor; $A =
I_0[T(q')f(q)]^2$, $q = 2 a_{f.c.c.} \sin \vartheta / \lambda$ and
$q' = q / a_{f.c.c.}$ are, respectively, the dimensionless and
the usual scattering vector length with $a_{f.c.c.}$ being the f.c.c.
bulk lattice constant; $N_{\mathcal X} (n)$ is the number of atoms
in the cluster, $u^{{\mathcal X},n}_{i,j}$ the distance between
the {\it i-th} and {\it j-th} atom, $a_{\mathcal X} (n)$ the strain
factor.
The total scattered intensity is computed as
\begin{equation}
I(q) = \sum_{\mathcal X} x_{\mathcal X} \sum_{n=1}^{S_{\mathcal
X}} f_{\mathcal X}(n) I_{{\mathcal X},n}(q), \label{eq:tot_i}
\end{equation}
where $S_{\mathcal X}$ denotes the size of the largest cluster of
type ${\mathcal X}$, $x_{\mathcal X}$ is the number fraction of
each structure type ($\sum_{\mathcal X} x_{\mathcal X}=1$), and
$f_{\mathcal X}(n)$ is the value of log-normal size distribution
(\ref{eq:nano_size_distr}). It is worth noting that both
intensities in (\ref{eq:debye}) and (\ref{eq:tot_i}) are actually
functions of the scattering angle $\vartheta$ being $q= 2
a_{f.c.c.} \lambda^{-1} \sin \vartheta$.
Experimental XRD intensity profiles are collected by counting, at
each scattering angle $\vartheta_n$, the number of scattered
photons giving the diffracted intensity signal $I_n$. For such
events data are affected by Poisson noise. Since the number of
photons scattered at each angle $\vartheta_n$ is large, the
Poisson-distributed noise can be approximated by a
Gaussian-distributed noise as required in section \ref{sec:HSVD}.

\begin{figure}
\begin{center}
   \includegraphics*[width=0.95\textwidth]{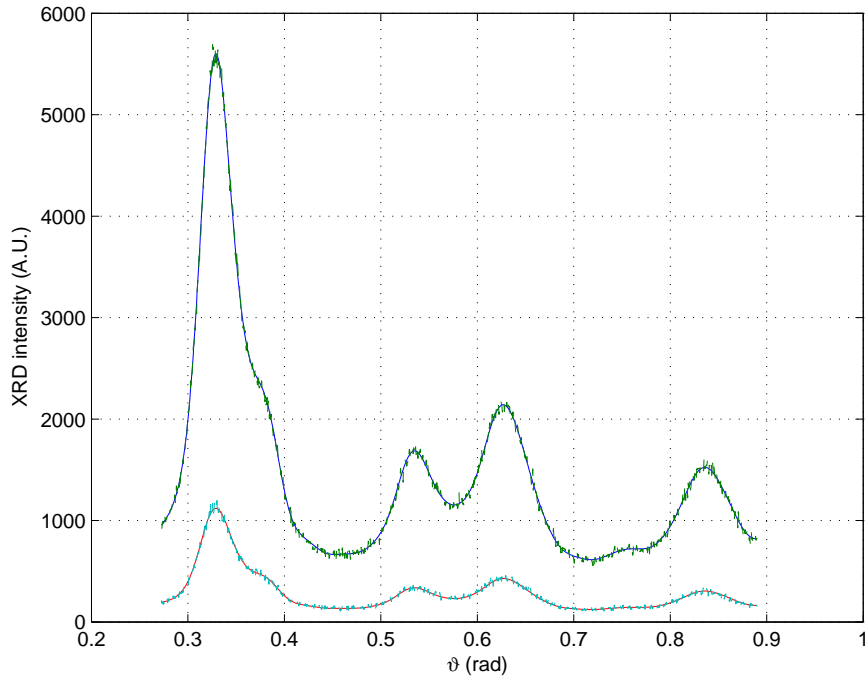}\\*
\end{center}
\caption{Synthetic XRD intensity profiles as a function of the
scattering angle. Table \ref{table_0} summarizes the parameters
used in eq.~(\ref{eq:debye}) to compute the diffraction
intensities. They are characteristics of Au samples. The
wavelength and bulk lattice constant have been set to
$\lambda=0.15418$ nm and $a_{f.c.c.}=0.40786$ nm, respectively.
From the upper to the lower profile, the NSR increases from $2$ to
$5\%$ (see fig.~\ref{fig_2} and text for details).}
\label{fig_1}
\end{figure}

\begin{figure}
\begin{center}
   \includegraphics*[width=0.95\textwidth]{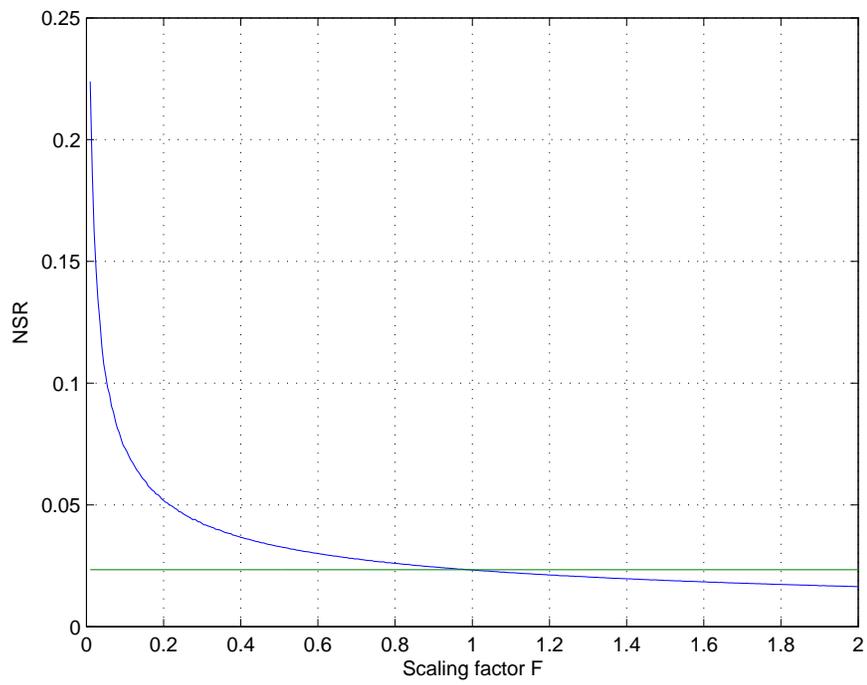}\\*
\end{center}
\caption{NSR as a function of the factor $F$. See text for
details.}
\label{fig_2}
\end{figure}

Noisy synthetic XRD intensity profiles were built by
generating Poisson distributed random profiles with intensity $I$
(\ref{eq:tot_i}) taken as the mean value of the Poisson process.
As a measure of the noise level, the noise-to-signal ratio (NSR)
was defined as follows:
\begin{equation}
{\rm NSR} = \frac{|| \sqrt{\mathcal{P} \left ( F I \right ) }
||}{|| \mathcal{P} \left ( F I \right )||}
                          \;,
\label{eq:nsr2}
\end{equation}
where $\mathcal{P}(I)$ denotes a Poisson process
with mean value $I$. Figure \ref{fig_1} displays XRD intensity
profiles with increasing NSRs. They were obtained by setting
$\lambda=0.15418$ nm and $a_{f.c.c.}=0.40786$ nm in
eq.~(\ref{eq:debye}). The set of parameters used to compute the
synthetic profiles are summarized in Table \ref{table_0}.
Different NSR values were obtained by scaling the scattered
intensity (\ref{eq:tot_i}) by a factor $F$. Figure \ref{fig_2}
shows the NSR of the synthetic profiles as a function of the
scaling factor F ranging from 0 to 2. This range contains the NSR
values usually measured in experimental profiles

We also considered real data in order to validate our method.
Three different samples were prepared with a resultant mean diameter of
2.0, 3.2 and 4.1 nm, respectively (as measured by TEM).
The size distributions were approximately characterized by the same width (FWHM
$\approx$ 1 nm) for all three systems.
Powder XRD studies were realized on the XRD beam line at the Brazilian
Synchrotron Light facility (LNLS--Campinas, Brazil) using 8.040 keV
photons at room temperature. For further details see Ref. \cite{zanchet00}.

\section{Numerical results}

Noisy synthetic XRD patterns were generated corresponding to
nanocrystalline samples of increasing size from 2 to 4 nm, and
Poisson-distributed noise with increasing NSR from 2\% to 10\%. The
HLSVD--PRO filter was then applied to the noisy synthetic XRD
signals in order to study its properties under controlled conditions. A key
step in the filtering procedure is the selection of the number $K$
of damped sinusoids characterizing the model function of the
HLSVD--PRO filter. Here, a possible approach is presented, which is based
on the following frequency selection criterion: the singular values
$\lambda_k$, $k=1,\dots,r$, are plotted vs. the corresponding
frequencies $f_k$ of the sinusoids in eq.~(\ref{eq:model}).
This choice facilitates a direct comparison of the results of the
proposed filter  with those obtained by other filters based on a
frequency approach. It was observed that, generally,
crystallographic XRD intensity signals show a clear transition
from a low-frequency region, characterized by high singular values $\lambda_k,$
to a high-frequency region with small singular values. The index
$K$ of the frequency $f_K$ corresponding to the transition,
provides the number of damped sinusoids to be used in the
HLSVD--PRO filter.

\begin{figure}
\begin{center}
   \includegraphics*[width=0.95\textwidth]{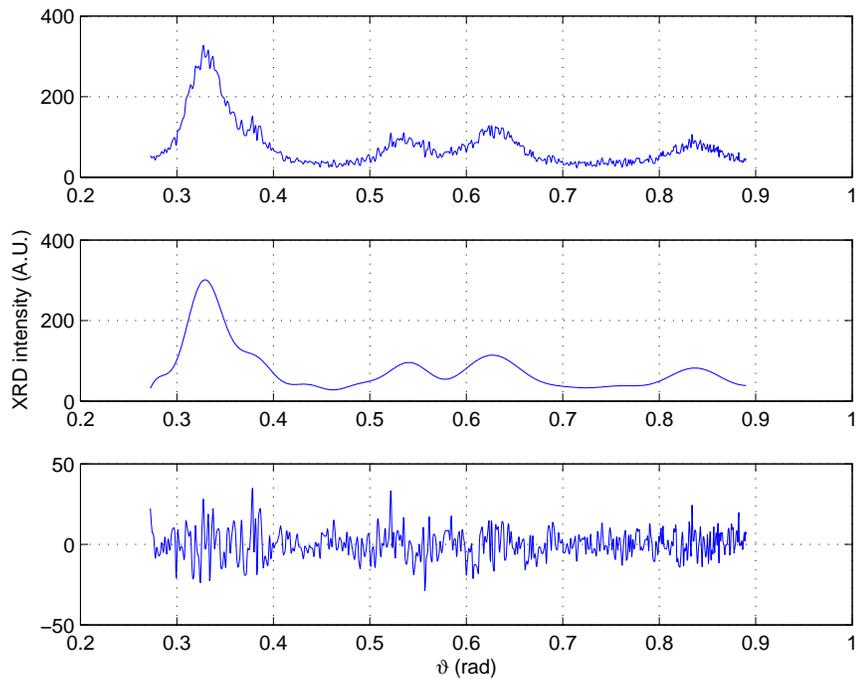}\\*
\end{center}
\caption{3 nm Au synthetic sample: (top) noisy (NSR=$10\%$)
synthetic XRD intensity profile as a function of the scattering
angle $\vartheta$; (middle) filtered XRD intensity profile. The
HLSVD--PRO filter removes signal components with frequency above
$f_K=35$ rad$^{-1}$ (see fig.~(\ref{fig_4})); (bottom)
difference between measured and filtered profiles.} \label{fig_3}
\end{figure}

\begin{figure}
\begin{center}
   \includegraphics*[width=0.95\textwidth]{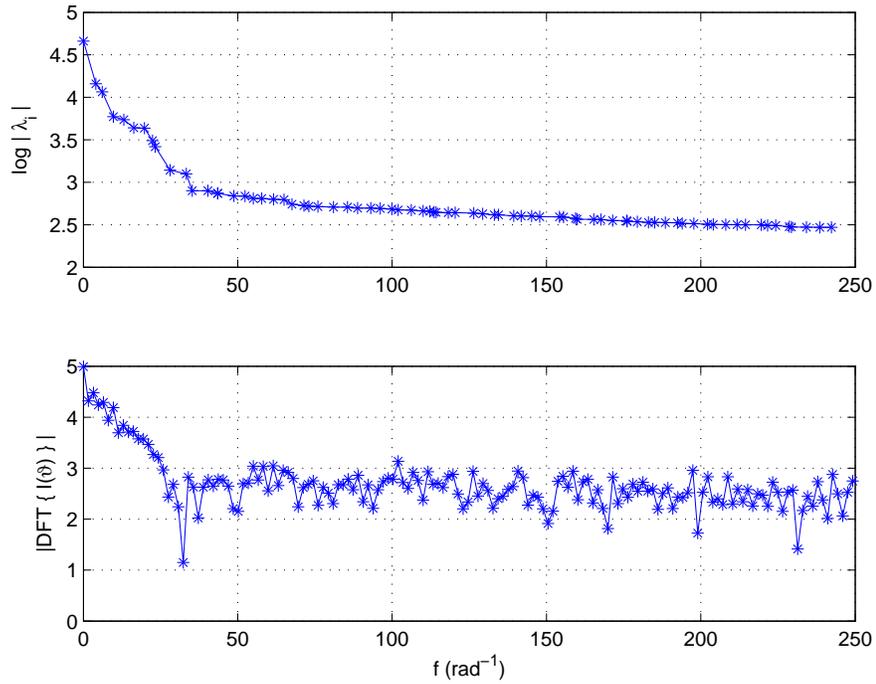}\\*
\end{center}
\caption{3 nm Au synthetic sample with NSR = $10\%$: (top)
amplitude of eigenvalues $\lambda_k$ vs. frequency $f_k$,
$k=1,\dots,r$. The frequency $f_K=35$ rad$^{-1}$ was used to
separate (filter) high frequency components due to noise; (bottom)
Portion of the DFT amplitude spectrum of the noisy synthetic XRD
intensity profile. Both plots refer to the XRD intensity profile
shown at the top of fig.~(\ref{fig_3}). } \label{fig_4}
\end{figure}

Figure \ref{fig_3}  displays an example of application of the
HLSVD--PRO filter. A noisy synthetic XRD intensity profile is shown at the 
top of the figure. 
It corresponds to X-ray scattering from a Au sample having a 3 nm size 
with a Poisson-like 
noise with NSR=10\%.   The filtered signal 
shown in the middle of the
figure was obtained by setting $K=9$ in the HLSVD--PRO filter. 
This value was estimated by visually inspecting the plot of the singular values
$\lambda_k$ vs the frequencies $f_k$ (see top of fig.~(\ref{fig_4})).
Specifically, a transition from high to small $\lambda_k$ was observed at
frequency $f_{K}=35$ rad$^{-1},$ which represents the $K^{th}$ frequency in the
set of the sorted frequencies starting from the smallest one. For a
comparison, the discrete Fourier transform (DFT) of the noisy
synthetic XRD signal is reported at the bottom of fig.~(\ref{fig_4}). 
Again, a phenomenon of transition from high to small
singular values occurs in the same region of the spectrum,
as observed at the top of figure. However, the transition frequency is much 
more
difficult to localize than in the HLSVD--PRO filter case. This
makes troublesome the application of DFT and WF filters to clean
noisy XRD data. It is worth noting that this
difference between the HLSVD--PRO and Fourier frequency based
filters is relevant when the filter is intended to be used during
interactive XRD data analysis. In this case the successful
application of a blackbox filter easy-to-use becomes crucial.

\begin{figure}
\begin{center}
   \includegraphics*[width=0.95\textwidth]{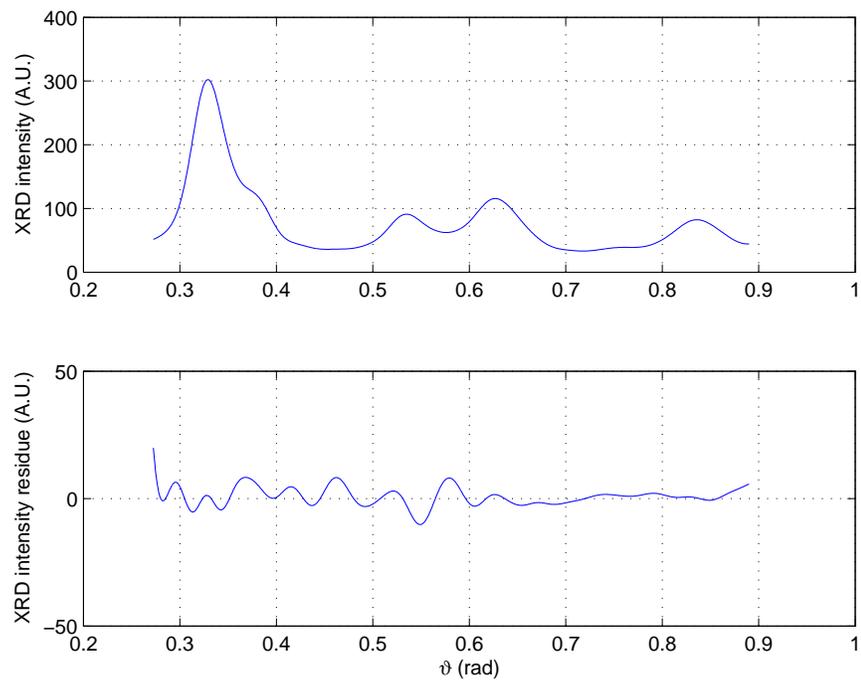}\\*
\end{center}
\caption{3 nm Au synthetic sample with NSR = $10\%$: (top)
noiseless synthetic XRD intensity profile as a function of the
scattering angle $\vartheta$; (bottom) difference between
noiseless and filtered (see middle plot of fig.~(\ref{fig_3}))
profiles.} \label{fig_5}
\end{figure}

\begin{figure}
\begin{center}
   \includegraphics*[width=0.95\textwidth]{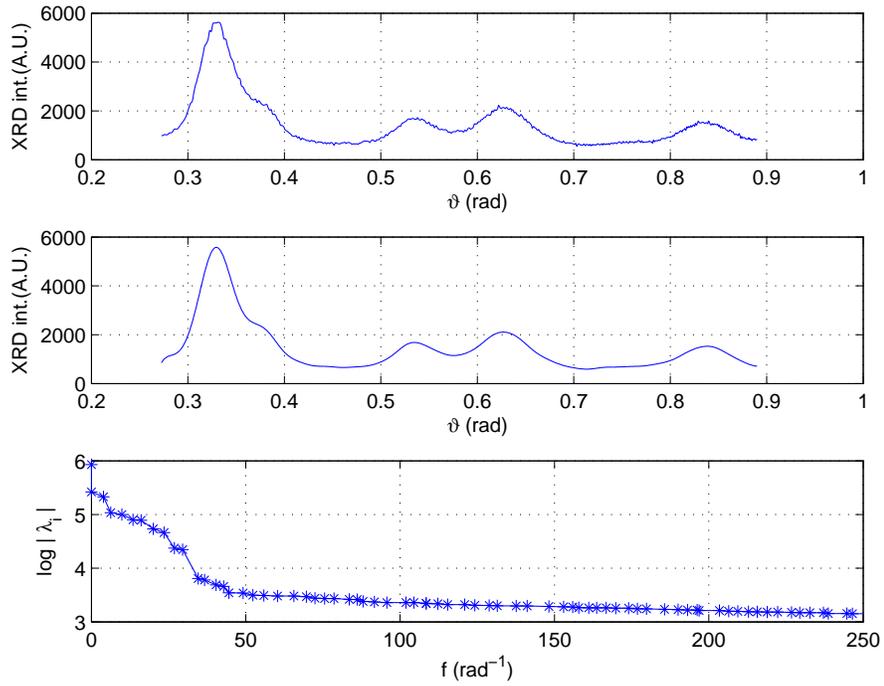}\\*
\end{center}
\caption{3.2 nm Au real sample: (top) noisy (NSR=2.3\%) XRD
intensity profile as a function of the scattering angle
$\vartheta$; (middle) filtered XRD intensity profile. The
HLSVD--PRO filter removes signal components with frequency above
$f_K=34$ rad$^{-1}$; (bottom) amplitudes of eigenvalues
$\lambda_k$ vs frequency $f_k$, $k=1,...,q$.
%
}\label{fig_6}
\end{figure}

Coming back to fig.~(\ref{fig_3}), the difference between the
values of noisy and filtered profiles
is shown at the bottom. To quantify the
performance of the filter, the filtered signal was compared
with the noiseless synthetic XRD signal (see
fig.~(\ref{fig_5})). This can be done only with synthetic
signals as experimental XRD data without noise are not available.
The filter performance was evaluated using the measure
\begin{equation}
{\mathcal E} = \frac{|| I_{exp} - I_{th} ||}{|| I_{fil} - I_{th} ||}
                          \;,
\label{eq:filter_perf}
\end{equation}
To give a statistical significance to these measures a Monte Carlo
experiment was carried out. More precisely, the HLSVD--PRO was applied
to 1000 noisy synthetic profiles generated by starting from the
same sample size, with different NSRs. For each filtered profile,
the measure ${\mathcal E}$ of filter performance was estimated
by calculating the mean value and the standard deviation. Table
\ref{table_1} summarizes the outcome of the aforementioned Monte Carlo 
experiments.
For each sample size and NSR, the mean and standard deviation are obtained 
using $1000$ synthetic XRD intensity profiles with different noise 
realizations having the same NSR.
The sensitivity to the number $K$ of sinusoids of the
HLSVD--PRO filter was also studied. This number was
slightly varied around the optimal $K$ value selected by using the
frequency criterion.
The performance results were compared in order to validate the choice of the 
optimal $K$ value. In particular, $K$ was increased and decreased by 2, 
as discussed
in Section I.
The results of such a comparison are summarized in Table \ref{table_2}
and show that the proposed frequency criterion provides the value of $K$
corresponding to the best performance of the HLSVD--PRO filter.
The filter was also applied to real XRD intensity profiles of
Au samples of size 2, 3.2 and 4.1 nm. Figure \ref{fig_6} shows at
top the profile of a 3.2 nm Au sample with $NSR=2.3\%$ computed as
$|| \sigma || / || I ||$, where $\sigma$ and $I$ are vectors with
the measured error and the intensity values, respectively. Since
in the case of XRD signals the noise follows the Poisson
distribution, $\sigma$ is given by $\sqrt{I}$. 
The result obtained by HLSVD--PRO is displayed in the middle of the figure. 
At bottom the plot of singular values is depicted vs. the frequency. 
Components with a frequency higher than $f_K=34$
rad$^{-1}$, due to noise, were removed. Denoising a real XRD
profile of 500 intensity data samples, as typical ones used in the
present study, requires about 11 seconds, using MatLab 7 on a
machine with a Intel Xeon 1.80 GHz processor and a 512 KB cache
size.

\section{Conclusions}

A filter based on the HLSVD--PRO method has been presented. It has
been applied to filter XRD patterns of nanocluster powders. The
filter performance has been studied on synthetic and real XRD
patterns with different NSRs. Results show that the proposed
filter is robust and computationally efficient. A further
advantage is its user-friendliness that makes it a useful blackbox
tool for the processing of XRD data.

\section{Acknowledgments}

We thank A. Cervellino, C. Giannini and A. Guagliardi for kindly
providing us with experimental XRD data.


\newpage

\begin{table}[ht!]
\caption{Values of parameters used in eq.~(\ref{eq:debye}) to
compute synthetic XRD intensity profiles. The wavelength and the
$f.c.c.$ bulk lattice constant were set to $\lambda =
0.15418$ nm and $a_{f.c.c.}=0.40786$ nm, respectively.}
\label{table_0} \centering

\vskip 2cm

\begin{tabular}{|c||c||c||c|}
\hline
 Parameter & ${\mathcal X}={\mathcal C}$ &
             ${\mathcal X}={\mathcal I}$ &
             ${\mathcal X}={\mathcal D}$\\
\hline
$\xi_{\mathcal X}$ & 5.0 & 5.0 & 5.0 \\
$s_{\mathcal X}$ & 0.3 & 0.3 & 0.3 \\
$n^0_{\mathcal X}$ & 4.0 & 4.0 & 6.0 \\
$\Omega_{\mathcal X}$ & 1.0 & 1.0 & 1.0 \\
$\Xi_{\mathcal X}$ & 1.0 & 1.0 & 1.0 \\
$w_{\mathcal X}$ & 0.5 & 0.5 & 0.5 \\
\hline
\end{tabular}
\end{table}

\newpage

\begin{table}[ht!]
\caption{Measure ${\mathcal E}$ of the filter performance (see eq.
\ref{eq:filter_perf}). For each sample size and NSR,
%
%
mean and standard deviation figures refer to Monte Carlo
experiments run on $1000$ synthetic XRD intensity profile with
different noise realizations having the same NSR.} \label{table_1}
\centering

\vskip 2cm

\begin{tabular}{|c||c||c||c|}
\hline
 & NSR=2\% & NSR=5\% & NSR=10\%\\
\hline 2 nm & $1.89 \pm 0.28$ & $2.34 \pm 0.18$ & $2.49 \pm 0.16$ \\
\hline 3 nm & $1.54 \pm 0.39$ & $1.87 \pm 0.16$ & $2.34 \pm 0.20$ \\
\hline 4 nm & $1.25 \pm 0.06$ & $1.56 \pm 0.12$ & $1.89 \pm 0.19$ \\
\hline
\end{tabular}
\end{table}

\newpage

\begin{table}[ht!]
\caption{Measure ${\mathcal E}$ (see eq. \ref{eq:filter_perf}) of
the filter performance as a function of the order $K$ of the
filter, namely the cutoff frequency $f_K$. The synthetic XRD
intensity data refer to different sample sizes and NSRs.The best
performance corresponds to the order $K$ reported in the middle
row of each NSR value. } \label{table_2} \centering

\vskip 2cm

\begin{tabular}{|c||c|c|c|c|}
\hline
         &            &       2 nm      &       3 nm      &      4 nm       \\
\hline \hline
         &     K-2    & $1.86 \pm 0.16$ & $1.25 \pm 0.09$ & $1.67 \pm 0.10$ \\
\cline{2-5}
NSR=10\% &     K=9    & $2.49 \pm 0.16$ & $2.34 \pm 0.20$ & $1.89 \pm 0.19$ \\
\cline{2-5}
         &     K+2    & $2.43 \pm 0.42$ & $2.28 \pm 0.21$ & $1.73 \pm 0.18$ \\
\hline \hline
         &     K-2    & $2.17 \pm 0.18$ & $1.81 \pm 0.16$ & $1.52 \pm 0.11$ \\
\cline{2-5}
NSR=5\%  &     K=11   & $2.34 \pm 0.18$ & $1.87 \pm 0.16$ & $1.56 \pm 0.12$ \\
\cline{2-5}
         &     K+2    & $2.22 \pm 0.28$ & $1.87 \pm 0.16$ & $1.48 \pm 0.09$ \\
\hline \hline
         &     K-2    & $1.80 \pm 0.21$ & $1.37 \pm 0.32$ & $1.13 \pm 0.14$ \\
\cline{2-5}
NSR=2\% &      K=15   & $1.89 \pm 0.28$ & $1.54 \pm 0.39$ & $1.25 \pm 0.06$ \\
\cline{2-5}
         &     K+2    & $1.86 \pm 0.18$ & $1.46 \pm 0.38$ & $1.18 \pm 0.09$ \\
\hline
\end{tabular}
\end{table}


\newpage

\newpage

\newpage

\newpage

\newpage

\newpage

\end{document}